\documentclass[11pt]{amsart}

 \usepackage{graphicx}
 \usepackage{amscd}
 \usepackage{mathtools}

 \def\a{\alpha}
 
 \def\be{\beta}
 
 \def\de{\delta}
 
 \def\e{\varepsilon}

 \def\ga{\gamma}
 \def\dga{{\dot{\gamma}}}

 \def\vr{\varphi}

 \def\om{\omega}
 \def\Om{\Omega}

 \def\re{{\mathbb R}}

 \def\then{\Longrightarrow}
 \def\ov{\overline}
 \def\Z{{\mathbb Z}}

 \def\D{{\mathbb D}}

 \def\H{{\mathbb H}}

 \def\bh{{\mathbf{h}}}

 \def\tH{{\widetilde H}}

 \def\tM{{\widetilde{M}}}
 \def\hM{{\widehat{M}}}

 \def\T{{\mathbb T}}

 \def\tq{{\tilde{q}}}
 
 \def\tq1{{\tilde{q}_1}}

 \def\oH{{\ov{H}}}
 \def\oL{{\ov{L}}}

 \def \lv{\left\vert}
 \def \rv{\right\vert}
 \def \lV{\left\Vert}
 \def \rV{\right\Vert}
 \def \ov{\overline}
 
 \def \then{\Longrightarrow}

 \DeclareMathOperator{\diam}{diam}

  \renewcommand{\proofname}{{\bf Proof:}}

 \theoremstyle{plain}

  \swapnumbers

 \newtheorem{Thm}{Theorem}[section]

 \newtheorem{Theorem}[Thm]{\bf Theorem}
 \newtheorem{Proposition}[Thm]{\bf Proposition}

 \theoremstyle{definition}

 \theoremstyle{remark}

 \newtheoremstyle{Cl}
  {5pt}
  {3pt}
  {\sl}
  {}
  {\it}
  {:}
  {.5em}
  {}

 \theoremstyle{Cl}

 \def\begincproof{
                  \renewcommand{\proofname}{\it Proof:}
                  \begin{proof}
                 }

 \def\endcproof{
                \renewcommand{\qedsymbol}{$\diamondsuit$}
                \end{proof}
                \renewcommand{\qedsymbol}{\openbox}
                \renewcommand{\proofname}{\bf Proof:}
               }

 \renewcommand{\proofname}{{\bf Proof:}}

 \title
 {Homogenization on manifolds}

 \author[G. Contreras]{Gonzalo Contreras}

\address{CIMAT \\
          A.P. 402, 36.000 \\
          Guanajuato. GTO \\
          M\'exico.}
\email{gonzalo@cimat.mx}
\thanks{Partially supported by CONACYT, Mexico, grant  178838.}

\begin{document}

\parskip +5pt

\begin{abstract}
We present a theorem by Contreras, Iturriaga and Siconolfi \cite{CIS} in which we
give a setting to generalize the homogenization of the Hamilton-Jacobi equation
 from tori to other manifolds.
\end{abstract}

\maketitle

A {\it homogenization problem} consists of  a Partial Differential Equation (PDE)
with a fast (oscillating) variable $\e$ and a slow variable. 
The homogenization result is that when the oscillating period $\e$ tends to
zero, there is a limit of the solutions $u^\e$ of the PDE to a solution of
an homogenized or ``averaged'' PDE.

An example of the homogenization result that we present here is the convergence
of the average distance in the universal cover of the torus $\T^2=\re^2/\Z^2$ to the distance 
in the  stable norm in $H_1(\T^2,\re)=\re^2$, when the diameter of the fundamental
domain $\e$ tends to zero (see fig. 1).
\begin{figure}[h]
  \hskip -3cm
\quad \resizebox*{12cm}{7cm}{\includegraphics{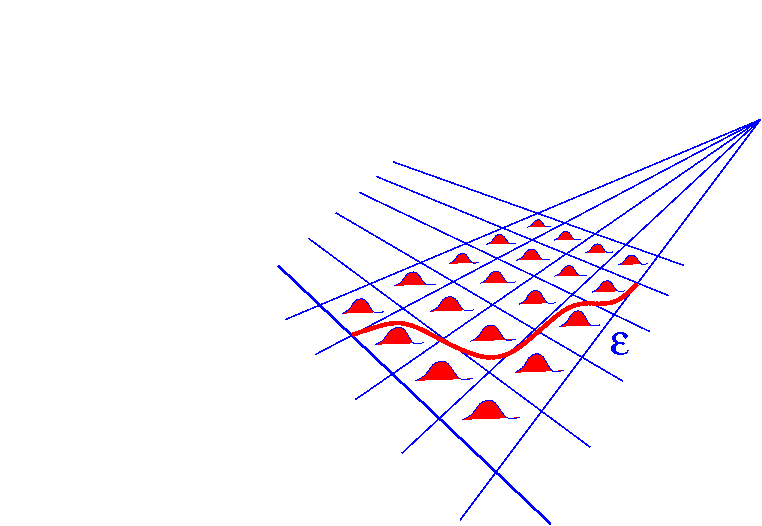}}
  \caption{Convergence to the stable norm.}
\label{mgeod}
\end{figure}  

\begin{figure}[h]
\hskip -3cm
\parbox{12cm}{
   \resizebox*{13cm}{7cm}{\includegraphics{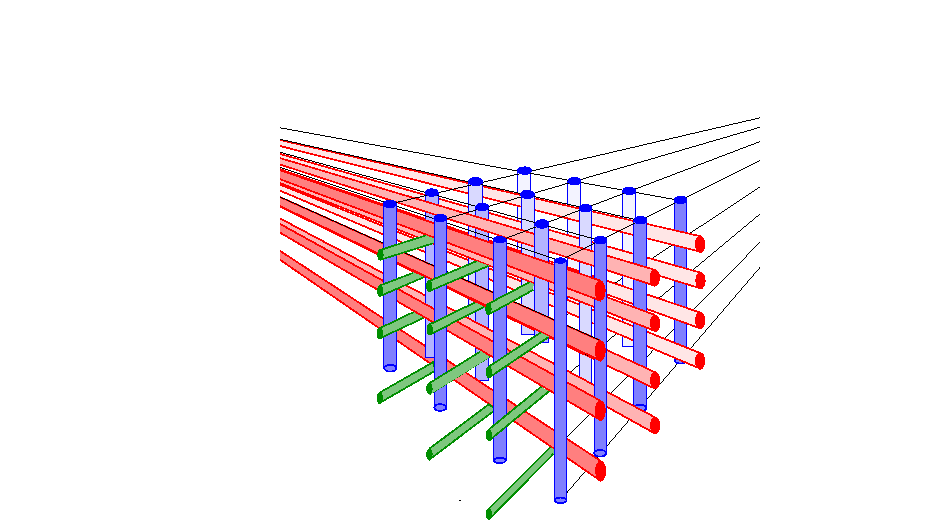} }
}
\parbox{2.5cm}{
   \resizebox*{2.5cm}{2.5cm}{\includegraphics{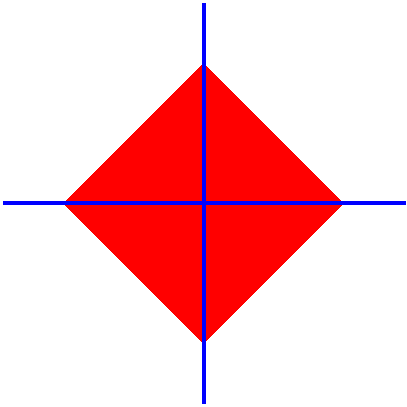}}
}
   \caption{Hedlund's example and its stable norm.}\label{fhedl}
\end{figure}

In higher dimensions the minimal geodesics may not converge. 
This is related to the flats of the stable norm as in
Hedlund's example \cite{Hedlund} in figure~\ref{fhedl}.

Hedlund's example is a 3-torus $\T^3=\re^3/\Z^3$, in which the  Riemannian metric 
is deformed in three disjoint tubes of different homological directions in which the
central closed geodesics are very short. In the example, minimal geodesics follow
the tubes with at most two jumps and the stable norm is
$$
\lV (x,y,z)\rV = |x|+|y|+|z|.
$$
In Hedlund's example the minimal geodesics do not converge as $\e\to 0$.
There is a convergence as ``holonomic measures'' to an invariant
measure supported on three periodic orbits on the tubes. The fact that 
there is no ergodic minimizing measure in a given homology
class implies that the stable norm is flat on that class.

 An important observation in this geodesic example of homogenization is that 
 the average minimal distance can be computed from the geodesics 
 of the stable norm, which are straight lines. One expects that the 
 homogenized or averaged problem is much simpler and computable 
 than the original problem. Another application of homogenization theory 
 is to obtain macroscopic laws from microscopic data.

 Homogenization theory has mostly been done in a periodic setting 
 (i.e. on the torus~$\T^n$) or in quasi-periodic tilings or random media 
 on $\re^n$.
 In the case of the Hamilton-Jacobi equation, the limiting objects are well 
 known and naturally defined on arbitrary manifolds: the effective Lagrangian
 is Mather's minimal action function $\beta$ and the effective (or homogenized)
 Hamiltonian is its dual $\beta^*$, also known as Ma\~n\'e's critical value.
  Nevertheless
 this homogenizations have only been made in~$\T^n$.
 
 We will show how to extend the homogenization result for the Hamilton-Jacobi equation from 
 the torus $\T^n$ to an arbitrary compact manifold. We hope that the setting presented here 
 can be applied to many other homogenization results.

 \section{Homogenization of the Hamilton-Jacobi equation.}

 Let $M$ be a compact manifold without boundary.
 A {\it Tonelli Lagrangian} is a $C^2$ function $L:TM\to\re$ satisfying:
 \begin{enumerate}
 \renewcommand\theenumi{\roman{enumi}}
 \item {\it Convexity:}  \quad $\frac{\partial L}{\partial v\,\partial v}(x,v)$ is positive definite $\forall (x,v)\in TM$.
 \item {\it Superlinearity:}\quad $\lim_{|v|\to+\infty}\frac{L(x,v)}{|v|}=+\infty$ uniformly on $x\in M$.
 \end{enumerate}

 Examples of Tonelli Lagrangians are 
 \begin{enumerate}
 \item {\it The kinetic energy: } $L(x,v)=\frac 12\,\lV v\rV_x$, which gives the geodesic flow and whose homogenization is equivalent to the examples given above.
 \item {The Mechanical Lagrangian:} $L(x,v)=\frac 12\,\lV v\rV_x -U(x)$ = kinetic energy - potential energy.
 This Lagrangian gives rise to Newton's law with force $F= -\nabla U(x)$.
 \end{enumerate}
 
 \medskip 
 
 The {\it action} of a smooth curve $\ga:[0,T]\to M$ is
 $$
 A_L(\ga)=\int_0^T L\big(\ga(t),\dot\ga(t)\big)\,dt.
 $$
 Critical points of $A_L$ satisfy the {\it Euler-Lagrange} equation
 \begin{equation}\label{EL}
 \frac{d}{dt} \frac{\partial L}{\partial v} = \frac{\partial L}{\partial x}.
 \end{equation}
 The Euler-Lagrange equation is a second order equation whose solutions 
 give rise to the {\it Lagrangian Flow:} $\vr_t:TM\to TM$, 
 $$
 \vr_t(x,v)=(\ga(t),\dga(t)),
 $$
 where $\ga$ is the solution of \eqref{EL} with initial conditions $(\ga(0),\dga(0))=(x,v)$.
 
 The convex dual of the Lagrangian is the Hamiltonian $H:T^*M\to\re$
 $$
 H(x,p) = \sup_{v\in T_xM} \big\{p(v) -L(x,v)\big\}.
 $$
 The {\it Legendre Transform} $L_v:TM\to T^*M$, $L_v(x,v)=\frac{\partial L}{\partial v}(x,v)$, converts the 
 Euler-Lagrange equation~\eqref{EL} into the Hamiltonian equations:
 $$
 \frac{d\,}{dt}\,L_v = L_x 
 \quad \xRightarrow{\;\;L_v\;\;} \quad
 \begin{cases}
 \dot x &= \phantom{-}H_p \\
 \dot p &= -H_x
 \end{cases}
 $$
 and conjugates the Lagrangian and Hamiltonian flows. 
 
 The {\it Hamilton-Jacobi} equation  
 \begin{equation}
 \partial_t u + H(x,\partial_x u)=0
 \end{equation}\label{HJ}
 encodes the minimal (Lagrangian) action cost.
 A solution $u:M\times \re_+\to\re$,
 to the Hamilton-Jacobi equation with initial condition
 $$
 u(x,0)=f(x)
 $$
 is given by the Lax formula
 $$
 u(x,t)=\inf\left\{f(\ga(0))+\int_0^t L(\ga,\dga)\;\Big\vert\; \ga\in C^1([0,t],M), \; \ga(t)=x\,\right\}.
 $$
 
 The characteristics of the Hamilton-Jacobi equation are {\it Tonelli minimizers} 
 i.e. minimizers of the action with fixed endpoints and fixed time interval.
 The value of the solution is the initial value + the action along these minimizers.
 Tangent vectors to the characteristics are related to $\partial_x u$ through the 
 Legendre Transform $L_v$:
 \begin{equation}\label{pxu}
 \partial_x u = L_v(\ga,\dga).
 \end{equation}
 
\begin{figure}[h]
\quad \resizebox*{8cm}{4.5cm}{\includegraphics{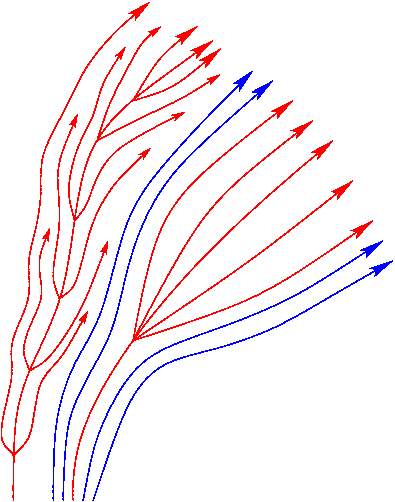}}
  \caption{Crossing of characteristics.}\label{crossingf}
\end{figure}

 Usually there are no global classical solutions of the Hamilton-Jacobi equation due to 
 crossing of characteristics as in figure~\ref{crossingf}.
 Indeed, from \eqref{pxu} at a crossing point there are various candidates for $\partial_xu$, 
 and hence $\partial_xu$ does not exist. 
   
 There are two popular types of weak solutions in PDEs:
 \begin{itemize}
 \item Weak solutions with {\it weakly differentiable} functions and 
         Sobolev Spaces are inspired  on the formula of integration by parts.
 \item The {\it viscosity solution} is inspired on the maximum principle for PDEs.
 \end{itemize}
 
 The first definition of viscosity solutions was made by L.C. Evans in 1980 \cite{Ev0}.
 Subsequently the definition and properties of the viscosity solutions of Hamilton-Jacobi equations
 were refined by Crandall, Evans and Lions in \cite{CEL}. 
 The existence and uniqueness of the {\it viscosity solution} of the initial value problem for the
 Hamilton-Jacobi equation was proved by Crandall and Lions in \cite{CrLi1}.
 
 A continuous function is a {\it viscosity solution} of
 $$
 \partial_t u + H(x,\partial_xu) =0
 $$
 if for every open set $U\subset M$ and any $\phi\in C^1(U\times\re_+,\re)$:
 \begin{itemize}
 \item if $u-\phi$ attains a local maximum at $(y_0,t_0)\in U\times\re_+$, then \\
 $\partial_t\phi(y_0,t_0)+H(y_0,\partial_x(y_0,t_0))\le 0$.
 \item if $u-\phi$ attains a local minimum at $(y_0,t_0)\in U\times\re_+$, then \\
 $\partial_t\phi(y_0,t_0)+H(y_0,\partial_x(y_0,t_0))\ge 0$.
  \end{itemize}

\begin{Theorem}[Lions, Papanicolaou, Varadhan \cite{LPV},
Evans \cite{Ev4}]\quad\label{LPV}

Let $H:\re^n\times\re^n\to\re$ be a $\Z^n$-periodic Tonelli Hamiltonian.
For $\e$ small let $f_\e:\re^n\to\re$ be Lipschitz. Consider the Cauchy
problem for the Hamilton-Jacobi equation
\begin{gather}
\partial_t u^\e + H\big(\tfrac x\e,\partial_xu^\e\big)=0,
\label{HJe}Ê\\
u^\e(x,0)=f_\e(x). \notag
\end{gather}
If $\lim_\e f_e = f$ uniformly then $\lim u^\e = u$ uniformly, 
where $u$ is the solution to
\begin{gather*}
\partial_t u +\oH(\partial_x u) =0,\\
u(x,0)=f(x).
\end{gather*}
The function $\oH:\re^n\to\re$, called the {\sl effective Hamiltonian} is 
convex, superlinear and is independent of the variable $x$.
\end{Theorem}

 The solutions to the homogenized problem can be easily written
 because the characteristics are straight lines and $p=\partial_xu$
 is constant along them
 \begin{equation*}
 \begin{cases}
 \dot p &= -\oH_x =0, \\
 \dot x &=\phantom{-} \oH_p =\text{constant}.
 \end{cases}
 \end{equation*}
 Thus
 $$
 u(y,t)=\min_{x\in\re^n}\big\{f(x)+t\,\oL\big(\tfrac{y-x}t\big)\big\},
 $$
 where
 $$
 \oL(x,v)=\max_{p\in\re^n}\big\{ p(v)-\oH(p)\big\}
 $$
 is the {\it Effective Lagrangian}.
 
 It turns out that the Effective Lagrangian $\oL=\be$ is Mather's minimal action function 
 $\be:H_1(\T^n,\re)\to\re$. The Effective Hamiltonian is related to Ma\~n\'e's critical value by
 $$
 \oH(P)=\a(P)=c(L-P), \qquad P\in\H^1(\T^n,\re),
 $$
 here $(L-P)(x,v):= L(x,v)-\om_x(v)$, 
 where $\om$ is a closed 1-form in the cohomology class~$P$.
 As such, it has several interpretations (see \cite{CILib}):
 \begin{enumerate}
 \renewcommand\theenumi{\roman{enumi}}
 \item\label{i} $\a$ is the convex dual of $\be$.
 \item $\a(P) = \inf\big\{\, k\in\re\;\vert\; \oint_\ga(L-P+k)\ge 0 \quad \forall\;
 \text{closed curve $\ga$ in $\T^n$}\,\big\}$.
 \item $\a(P)=\inf\big\{\, k\in \re\;\vert\; \Phi_k>-\infty\,\big\}$,\quad 
 where $\Phi_k:M\times M\to \re$ is 
 $$
 \Phi_k(x,y):=\inf\big\{\,\textstyle{\oint}_\ga(L-P+k)\;\vert\; 
 \ga \text{ curve in $\T^n$ from $x$ to $y$ }\big\},
 $$
 i.e. the minimal action with free time interval.\footnote{The function $\Phi_k$ is
 called Ma\~n\'e's action potential.}
 \item $\a(P)=-\inf\big\{\textstyle\int (L-P)\,d\mu 
 \;\vert\; \mu\text{ is an invariant measure for }L\;\big\}$.
 \item $\a(P)$ is the energy level  containing the support of the invariant measures $\mu$ which minimize $\int (L-P)\,d\mu$.
 \item $\a(P)= \displaystyle 
 \min_{u\in C^1(\T^n,\re)}\; \max_{x\in \T^n} \; H(x,P+d_xu)$.
 \item $\a(P)$ is the minimum of the energy levels which contain a Lagrangian graph in $T^*\T^n$ 
 with cohomology class $P$.
 \item\label{wKAM} From Fathi's weak KAM theory, $\a(P)$ is the unique constant for which there are global 
 viscosity solutions of the Hamilton-Jacobi equation
 $$
 H(x,P+d_xv)=\a(P), \qquad x\in \T^n.
 $$
 \end{enumerate}
 
 We explain briefly why Theorem~\ref{LPV} and \eqref{wKAM} 
 imply that the Effective Hamiltonian $\oH$ is
 Mather's alpha function $\a$. 
 Consider the case of affine initial conditions. The problem
 \begin{equation}\label{paffine}
 \left\{
 \begin{aligned}
 f(x)= u(x,0)&=a+P\cdot x\;\\
 \partial_tu +\a(\partial_xu)&=0
 \end{aligned}
 \right\}
 \end{equation}
 has solution 
 $$
 u(x,t)=a+P\cdot x-\a(P) t.
 $$
 Let $v:\T^n\times\re_+\to\re$ be a $\Z^n$-periodic solution to the {\it ``cell problem'':}
 $$
 H(x,P+d_xv)=\a(P), \qquad v:\T^n\times\re_+\to\re.
 $$
 Let 
 \begin{align*}
 u^\e(x,t)&:= u(x,t)+\e\,v\big(\tfrac x\e\big),\\
 F_\e(x) &:= u^\e(x,0)= f(x)+\e\,v\big(\tfrac x\e\big).
 \end{align*}
 Then $u^\e$ solves
 \begin{equation*}
 \left\{\;
 \begin{aligned}
 &\partial_tu^\e+H\big(\tfrac x\e,\partial_xu^\e\big)
 =-\a(P)+H\big(\tfrac x\e,P+\partial_y v\big(\tfrac x\e\big)\big) = 0,
 \\
 &u^\e(x,0)=f_\e(x).
 \end{aligned}
 \right.
 \end{equation*}
 Also we have that $f_\e\to f$ and $u^\e \to u$ uniformly and by \eqref{paffine} $u$
 satisfies a Hamilton-Jacobi equation with Hamiltonian $\a$. Therefore
 Theorem~\ref{LPV} implies that $\oH(P)=\a(P)$.
 
 \bigskip
 
 \subsection{The Problems}\quad 
 
 \bigskip
 
 The generalization of Theorem~\ref{LPV} to other manifolds has three problems:
 \begin{itemize}
 \item[1.]\label{p1} It is not clear how to choose the generalization of $\frac x\e$.
 \item[2a.]\label{p2a} Equation \eqref{HJe} is the Hamilton-Jacobi equation for the Hamiltonian
 \linebreak
 $H_\e(x,p):= H(\tfrac x\e,p)$, where $p$ {\it ``remains the same''}.  It is not clear how to do 
 it in non-parallelizable manifolds where the parallel transport depends on the path.
 \item[2b.]\label{p2b} The effective Hamiltonian $\oH(P)$ ``does not depend on $x$''. 
 This is another version of the same problem 2a.
 \item[3.] The candidate for effective Hamiltonian is Mather's $\a$ function $\a:H^1(M,\re)\to\re$.
 But in general $\dim H^1(M,\re) \ne \dim M$, i.e. the limit PDE would be in a space with different 
 dimension, the differential structure would be destroyed.
\end{itemize}
 
 \medskip
 
 In fact, the Hamilton-Jacobi equation is an encoding of a variational principle 
 (the minimal cost function) that will be stable under the change of space. 
 
 The torus $M=\T^n$ has many coincidences that allow to formulate Theorem~\cite{LPV}:
 \begin{enumerate}
 \item Its universal cover satisfies
 $$
 \widetilde{\T^n} =\re^n=H_1(\T^n,\re)=H^1(\T^n,\re).
 $$
 The effective Hamiltonian $\oH=\a:H^1(\T^n,\re)=\re^n\to\re$ and the
 effective Lagrangian $\oL=\be:H_1(M,\re)\to\re$ are defined in the 
 same space as the original periodic Hamiltonian.
 Thus the original PDE and the limit equation are in the same space.
 
 \item The cotangent bundle is trivial: $\text{T}^*\T^n=\T^n\times\re^n$ 
 and the parallel transport does not depend on the path. 
 Thus we can talk of a Hamiltonian that does not depend on $x$ and the 
 Hamilton-Jacobi equation for the effective Hamiltonian $\partial_tu +\oH(\partial_x u) =0$
 makes sense.
 
 \end{enumerate}

 \bigskip
 
 \subsection{The solution.}\quad
 
 \bigskip
 
 {\bf Problem 3.} We start with the solution to problem 3: the space for the family of PDEs.
 Let $M$ be a compact manifold without boundary and $H:T^*M\to \re$ a Tonelli Hamiltonian.
 Consider the Hurewicz homomorphism $\bh:\pi_1(M)\to H_1(M,\re)$ which sends the homotopy
 class of a curve to its homology class with real coefficients. The maximal free abelian cover $\tM$
 is the covering map $\tM\to M$ with group of Deck transformations 
 $$
 \text{Deck}(\tM)=\Z^k=\text{Im}(\bh)\subset H_1(M,\re),
 $$
 where $k=\dim H_1(M,\re)$ and $\pi_1(\tM)=\ker\bh$.

 \bigskip
 
{\bf Problem  1.[$x\mapsto\tfrac x\e$]} 
Let $d$ be the metric induced on $\tM$ by the lift of the Riemannian metric on $M$.
 For problem 1   we use the metric spaces $M_\e:=(\tM,\e d)$.
 
\begin{figure}[h]
  \quad \resizebox*{6cm}{4cm}{\includegraphics{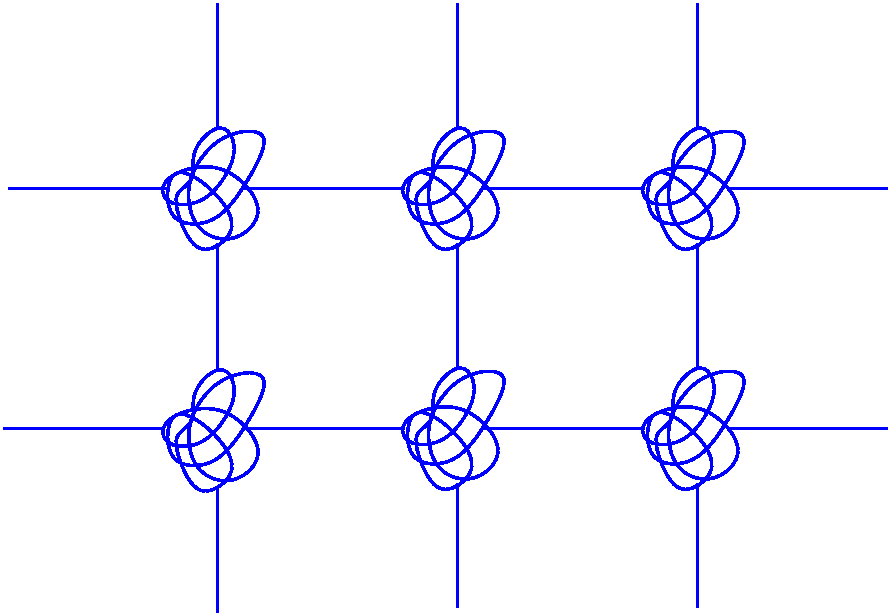}}
\caption{The structure of $\tM$.}\label{stru1}
\end{figure}

 The maximal free abelian cover $\tM$ has the structure of $\Z^k$, i.e. it is 
 (perhaps a complicated)  fundamental  domain
 which is repeated as the points in $\Z^k$, as in figure~\ref{stru1}.
The space $\tM_\e$ has a ``large scale structure'' as 
$\e\Z^k\hookrightarrow \re^k=H_1(M,\re)$.

We think of $M_\e\xrightarrow{\,\;\e\;\,} H_1(M,\re)$ as of $\e\Z^k\longrightarrow \re^k$.
For example: ``linear maps on $\tM$'' shall correspond to  ``integrals of closed 1-forms''.
Our solutions of the $\e$-oscillation Hamilton-Jacobi equation will be uniformly Lipschitz 
on $\tM_\e$, i.e. $\e K$-Lipschitz on $M_\e$. So that a solution $U^\e$ on $M_\e$ will 
define a function $r
v^\e$ on $\e \Z^k$ which is $K$-Lipschitz. By an Arzel\'a-Ascoli argument
we will obtain a convergence $v^\e\to v$ on $\re^k=H_1(M,\re)$.

\begin{figure}[h]
  \quad \resizebox*{5cm}{5cm}{\includegraphics{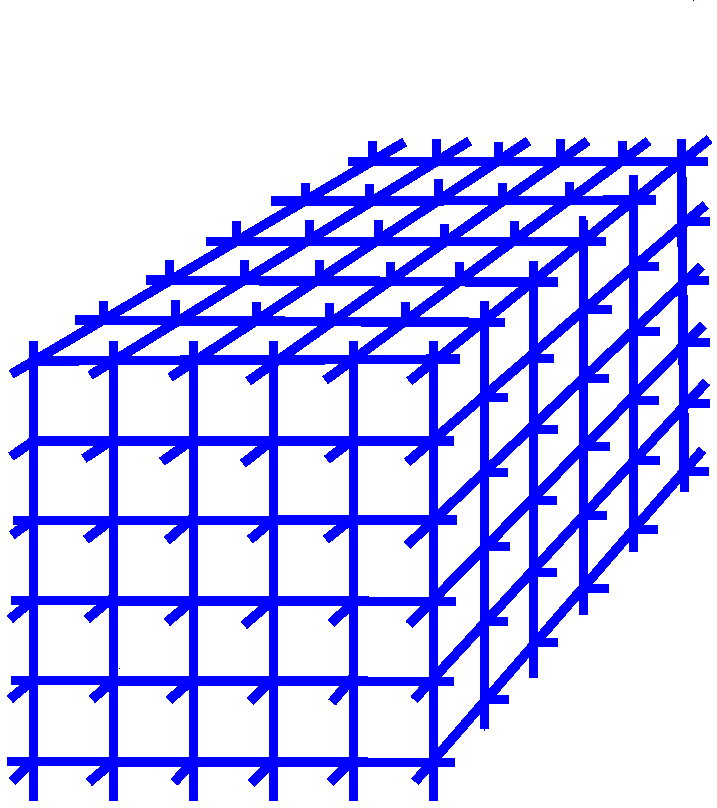}}
  \caption{Example of a free abelian cover of a surface $M=\T^2\#\T^2$ with group 
                  of Deck transformations $\Z^3$. It is not the maximal free abelian cover of $M$, 
                  because $\dim H_1(M,\re)=4$. The limit space $\lim_\e M_\e=\re^3$ has
                   higher dimension than $M$.}
  \label{3d}                 
\end{figure}

\bigskip

{\bf Problem 2. [$\oH$ independent of $x$] }
The solution to problem 2 consists on transforming the equation to an equivalent PDE.
In the case of $\re^n$ as in problem~\eqref{HJe} define $v^\e:\re^n\times\re^+\to\re$ by
$$
u^\e(x,t)=:v^\e(\tfrac xe,t).
$$
From~\eqref{HJe} we obtain that $v^\e$ is a solution to the problem
\begin{gather}
\partial_t v^\e+H\big(y,\tfrac 1\e\,\partial_y v^\e\big)=0, 
\label{HJev} \\
v^\e(y,0)=f_\e(\e y).
\label{vin}
\end{gather}
Now equation~\eqref{HJev} makes sense on any manifold.
Equation~\eqref{vin} will make sense with the
 following definition of convergence of spaces.

\bigskip

\subsection{Convergence of spaces.}\quad

\bigskip

This is inspired in Gromov's Hausdorff  convergence but it is 
made ad hoc for our homogenization problem. We will only need 
quasi-isometries because since we are doing analysis, just the 
equivalence class of the norms matter.

Let $(M,d)$, $(M_n,d_n)$  be metric spaces and $F_n:(M_n,d_n)\to (M,d)$ 
a continuous function. We say that $\lim_n(M_n,d_n,F_n)=(M,d)$ if
\begin{enumerate}
\renewcommand\theenumi{\alph{enumi}}
\item
There are $B, A_n>0$, with $\lim_n A_n=0$ such that 
$$
\forall x, y \in M_n: \qquad
B^{-1}\, d_n(x,y) - A_n \le d\big(F_n(x),F_n(y)\big)
\le B d_n(x,y).
$$
\label{ma}
\item For all $y\in M$ and $n$ there are $x_n\in M_n$ with $\lim_n x_n =y$.
\label{mb}
\end{enumerate}

Observe that \eqref{mb} is a kind of surjectivity condition. And \eqref{ma} implies 
that  
$$
\forall y\in M:
\qquad
\diam F_n^{-1}\{y\}\le B \, A_n \xrightarrow{\;\;n\;\;} 0,
$$
a kind of injectivity condition.

If $\lim_n(M_n,d_n,F_n)=(M,d)$, and $f_n(M_n,d_n)\to\re$, $F(M,d)\to\re$ are continuous,
we say that $\lim_n f_n=f$ {\it uniformly on compact sets} if for every compact set $K\subset M$
$$
\lim_n\sup_{x\in F_n^{-1}(K)} 
\lv f_n(x) - f(F_n(x))\rv 
=0.
$$
And we say that the family $\{f_n\}$ is equicontinuous if for every $\e>0$ there is $\de>0$
such that 
$$
\forall n: \qquad x,y\in M_n, \quad
d_n(x,y)<\de \quad \then \quad 
\lv f_n(x)-f_n(y)\rv <\e.
$$

\bigskip

Fix a basis $c_1,\ldots c_k$ for $H^1(M,\re)$.
Fix closed 1-forms $\om_i$ on $M$ such that $c_i=[\om_i]$.
Define $G:\tM\to H_1(M,\re)= H^1(M,\re)^*$ by
$$
G(x)\cdot c_i = \oint_{x_0}^x \tilde{\om_i}\;,
$$
where $\tilde\om_i$ is the pullback of $\om_i$ on $\tM$.
Let $F_\e:(M_\e,d_\e)\to H_1(M,\re)$ be $F(x):= \e\, G(x)$.
    
\begin{Proposition} 
$\lim_{\e\to 0}(\tM,\e d, F_\e)= H_1(M,\re)$
\end{Proposition}    
    
In the  homogenized or averaged problem we will have that    
the (limit) positions are in the configuration space $H_1(M,\re)$
and the momenta $p$ and differentials $\partial_xu$ are in the
dual of the configuration space $H_1^* = H^1(M,\re)$.

This explains why the effective Lagrangian $\oL=\be:H_1(M,\re)\to\re$
is defined in the {\it homology} group  $H_1(M,\re)$ but the
effective Hamiltonian is defined in the {\it cohomology} group
$H^1(M,\re)$.

\begin{Theorem}[Contreras, Iturriaga, Siconolfi \cite{CIS}]\label{T1}\quad

Let $M$ be a closed Riemannian manifold.
Let $H:T^*M\to\re$ be a Tonelli Hamiltonian and
$f_\e:(M_\e,d_\e)\to \re$ continuous functions such that
$\lim_\e f_\e= f$ uniformly, with $f:H_1(M,\re)\to\re$ Lipschitz. 

Let $\tH$ be the lift of $H$ to $\tM$ and let $v^\e$ be the
solution to the problem
\begin{gather*}
\partial_t v^\e+\tH\big(y,\tfrac 1\e \partial_y v^\e\big) =0, \\
v^\e(y,0)=f_\e(y).
\end{gather*}
Then the family $v^\e:\tM_\e\times]0,+\infty[\to\re$ is 
equicontinuous and
 $$
 \lim_{\e\to 0}v^\e=u:H_1(M,\re)\to\re
 $$
uniformly on compact sets of $H_1(M,\re)\times]T_0,+\infty[$,
for any $T_0>0$,
where $u$ is the solution to
\begin{gather*}
\partial_t u +\oH(\partial_x u) = 0,\\
u(x,0)=f(x);
\end{gather*}
and $\oH:H^1(M,\re)\to\re$ is $\oH=\a$ Mather's alpha function.
\end{Theorem}

 \subsection{Subcovers}\quad
 
\begin{figure}[h]
  \resizebox*{2cm}{5cm}{\includegraphics{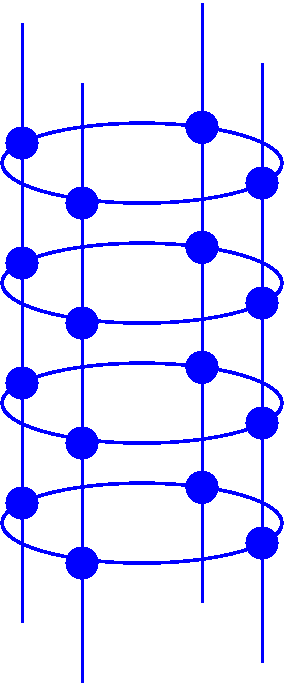}}
 \caption{The limit process kills the torsion: $\e\,(\Z_4\oplus \Z) \to \re$.}
 \label{torsion2}
 \end{figure}

  On abelian covers $\hM$ with Deck transformation group $\D$ of the form
 \linebreak
 $\D=\Z^k\oplus  \Z_{a_1} \oplus\cdots\oplus\Z_{a_p}$
 the limit $\lim_\e(\hM,\e d)$ will kill the torsion
 $\Z_{a_1} \oplus\cdots\oplus\Z_{a_p}$, as in figure~\ref{torsion2}.
 Thus we may restrict to free abelian covers with
 group of Deck transformation without torsion $\D=\Z^k$.
 These are sub covers of $\tM$.
 
 Using equivariance properties of the Hamilton-Jacobi equation, we
 obtain as a corollary of Theorem~\ref{T1} a similar result for other
 free abelian covers.
 
 \begin{figure}[h]
  \quad \resizebox*{6cm}{5cm}{\includegraphics{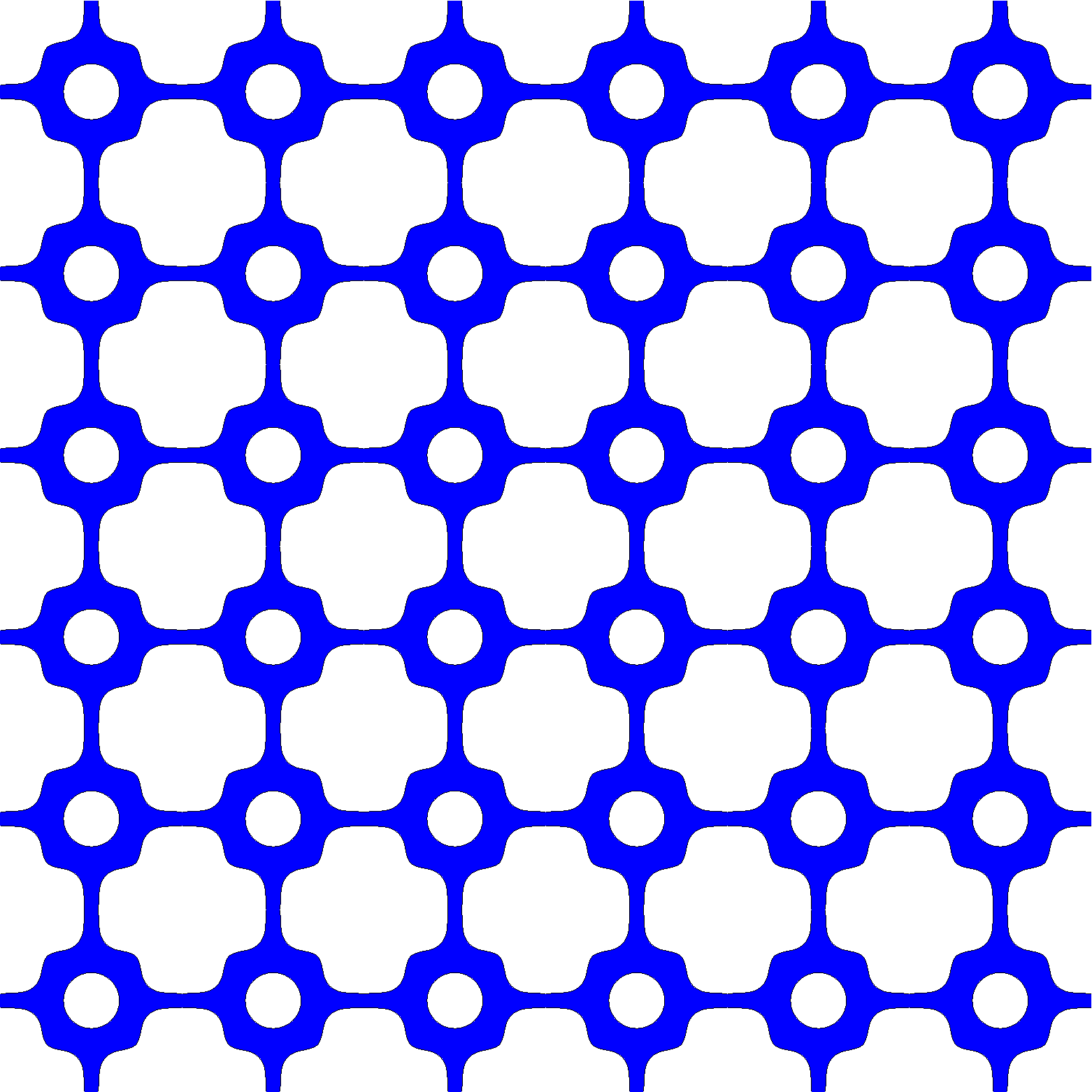}}
  \caption{$\Z^2$-cover of the surface $M_3$ of genus 3. In this case 
  $H_1(M_3,\re)=\re^3$.}
  \label{cover-talk}
\end{figure}

    \bigskip
    
    \subsection{Speculations.}\quad
    
    \bigskip
    
    There are generalizations of Aubry-Mather theory which can be interpreted 
    as a homogenization besides $\T^n$ or $\Z^n$ 
    and should give results in the setting presented above.
    On a generalization originated by Moser~\cite{MoserHD}, Caffarelli, de la Llave 
    and Valdinocci extend Aubry-Mather theory to higher dimensions 
    on very general manifolds, see \cite[remark 2.6]{LLV1}, \cite{calla1}, \cite{Caffa1}.
    There is also an extension by Candel and de la Llave \cite{CLL1} of the Aubry-Mather
    theory in statistical mechanics to configuration sets more general than $\Z^n$.
    Viterbo's symplectic homogenization \cite{ViterboSH} has also been extended to general
    manifolds by Monzner, Vichery and Zapolsky~\cite{MVZ}.
    
    Most of the homogenization theory is made only for the torus $\T^n$.
    Some PDE's techniques go through this setting despite the destruction of the 
    differential structure in the limit. 
    For example in the homogenization of the Hamilton-Jacobi equation,
    Evans perturbed test function method goes through to give a proof of the
    same result.

    The translation of homogenization results to manifolds can give interesting geometric objects.
    We have the following examples:
    \begin{itemize}
    \item The homogenization of the geodesic flow gives the stable norm.
    
    The stable norm was used by Burago and Ivanov  in their proof of
    the Hopf conjecture~\cite{BuIv}. Bangert \cite[Th. 6.1]{bangert2} proves that 
    a metric on  $\T^2$ whose stable norm is euclidean is the flat metric on  $\T^2$.
    Osuna \cite{Osuna2} proves that if $\T^n$ has the 1-dimensional and $(n-1)$-dimensional 
    stable norms Euclidean then the metric is flat.
    
    \medskip
     
    \item The homogenization of the Hamilton-Jacobi equation gives Mather's alpha function or
    Ma\~n\'e's critical value as the effective Hamiltonian. 
    
    In this case the limiting object $\oH(P)$ 
    was known  independently of homogenization and had many interesting characterizations 
    besides homogenization: variational, ergodic, geometric, symplectic as in \eqref{i}--\eqref{wKAM}.
    \end{itemize}

    Another example of a possible result is the homogenization of the Riemannian Laplacian.
    Let $M$ be a closed manifold and $\Om\subset H_1(M,\re)$ a domain. Let $f:\partial \Om\to\re$
    and $F:\Om\to\re$ be continuous functions. Choose a basis $[\om_i]$ for $H^1(M,\re)$ and let
    $G_\e:\tM_\e\to H_1(M,\re)$ be
    $$
    G_\e(x)\cdot[\om_i]=\e \oint_{x_0}^x\om_i.
    $$
    \begin{figure}[h]
     \resizebox{7cm}{7cm}{\includegraphics{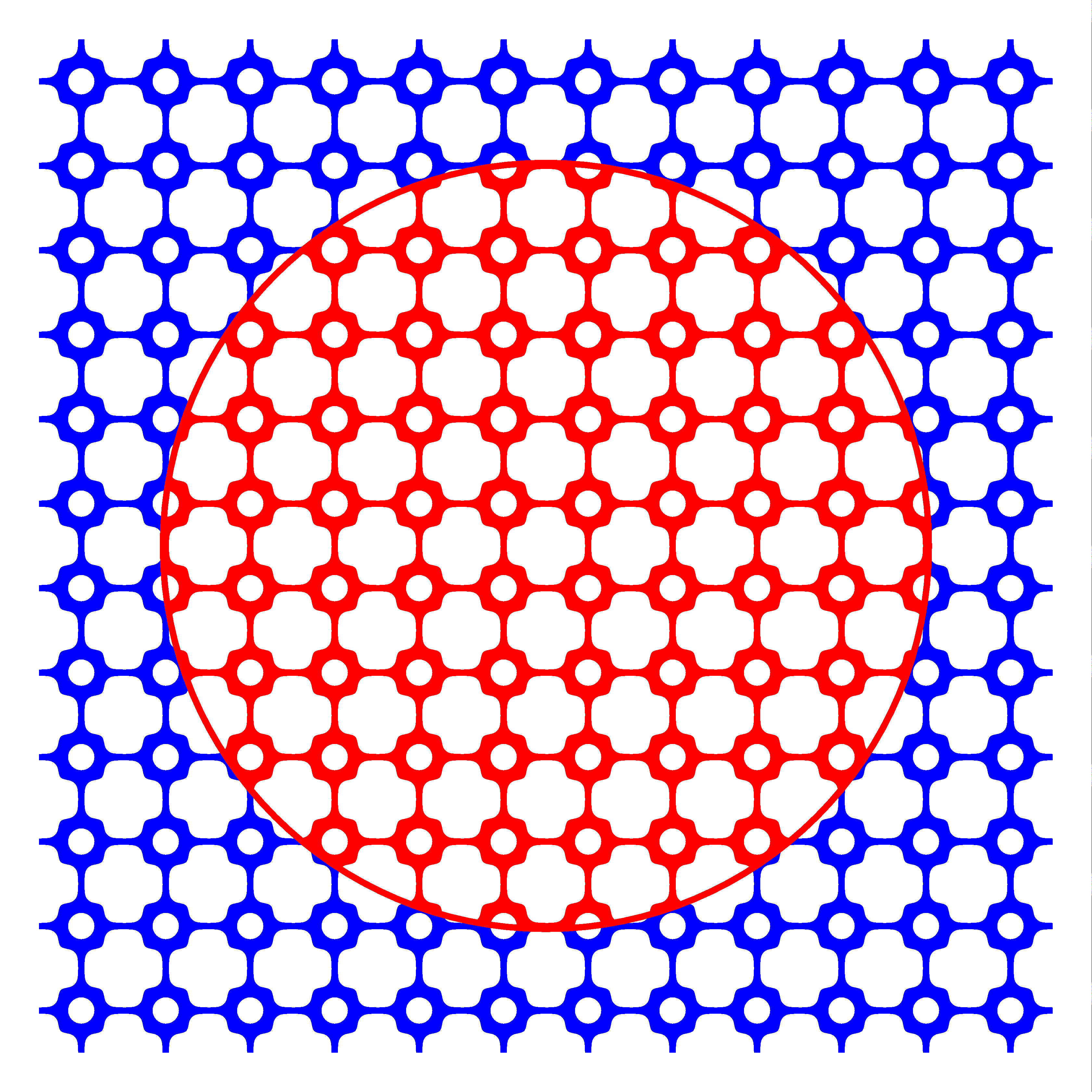}}
     \caption{Homogenization of the Riemannian Laplacian.}
    \label{lapla4}
    \end{figure}  
    Let $v_\e$ be the solution to the problem
    \begin{alignat*}{3}
    \Delta v_\e &= F\circ G_\e  \qquad &&\text{on } \quad G_\e^{-1}(\Om),Ê\\
    v_\e &= f\circ G_\e  &&\text{on } \quad \partial G_\e^{-1}(\Om).
    \end{alignat*}
    Prove that $v_\e \to u$ where
    \begin{alignat*}{2}
    \sum_{ij} A_{ij} \,\frac{\partial^2 u}{\partial x_i\partial x_j} &= F \qquad &&\text{on }\quad \Om, \\
    u &= f &&\text{on } \quad \partial \Om.
    \end{alignat*}
    In this homogenized Laplacian we should have that
    $$
    A_{ij} =\int_M \langle \eta_i(x),\eta_j(x)\rangle 
    $$
    where $\langle\cdot,\cdot\rangle$ is the induced inner product in $T^*M$ and
    $\eta_i$ is the harmonic 1-form in the class $[\om_i]$.
    
    \pagebreak
    
    Other questions can be:
    \begin{itemize}
    \item Homogenization of the eigenvalue problem for the Riemannian Laplacian.
    \item Probabilistic proofs of the homogenization of the Laplacian.
    \item Homogenization of the discretization of the Laplacian on graphs.
    \item Does it always give the same effective Laplacian? \\
            Also for the wave and heat equations? 
    \item What about quasi-periodic arrays of manifolds?
    \item What about non-abelian covers?        
    \end{itemize}
    
    For non-abelian covers we have some work in progress with Alfonso Sorrentino.
    The Gromov-Hausdorff tangent cone of the covering \cite{BBI} should give the effective space.


\begin{thebibliography}{10}

\bibitem{bangert2}
Victor Bangert, \emph{{Geodesic rays, Busemann functions and monotone twist
  maps.}}, Calc. Var. Partial Differ. Equ. \textbf{2} (1994), no.~1, 49--63
  (English).

\bibitem{BuIv}
D.~Burago and S.~Ivanov, \emph{Riemannian tori without conjugate points are
  flat}, Geom. Funct. Anal. \textbf{4} (1994), no.~3, 259--269.

\bibitem{BBI}
Dmitri Burago, Yuri Burago, and Sergei Ivanov, \emph{A course in metric
  geometry}, Graduate Studies in Mathematics, vol.~33, American Mathematical
  Society, Providence, RI, 2001.

\bibitem{Caffa1}
Luis Caffarelli, \emph{A homogenization method for non variational problems},
  Current developments in mathematics, 2004, Int. Press, Somerville, MA, 2006,
  pp.~73--93.

\bibitem{calla1}
Luis~A. Caffarelli and Rafael de~la Llave, \emph{Planelike minimizers in
  periodic media}, Comm. Pure Appl. Math. \textbf{54} (2001), no.~12,
  1403--1441.

\bibitem{CLL1}
A.~Candel and R.~de~la Llave, \emph{On the {A}ubry-{M}ather theory in
  statistical mechanics}, Comm. Math. Phys. \textbf{192} (1998), no.~3,
  649--669.

\bibitem{CILib}
Gonzalo Contreras and Renato Iturriaga, \emph{{G}lobal {M}inimizers of
  {A}utonomous {L}agrangians}, $22^{\text{o}}$ Coloquio Bras. Mat., IMPA, Rio
  de Janeiro, 1999.

\bibitem{CIS}
Gonzalo Contreras, Renato Iturriaga, and Antonio Siconolfi,
  \emph{Homogenization on arbitrary manifolds},  (2014), to appear in Calculus
  of Variations and P.D.E.

\bibitem{CEL}
M.~G. Crandall, L.~C. Evans, and P.-L. Lions, \emph{Some properties of
  viscosity solutions of {H}amilton-{J}acobi equations}, Trans. Amer. Math.
  Soc. \textbf{282} (1984), no.~2, 487--502.

\bibitem{CrLi1}
Michael~G. Crandall and Pierre-Louis Lions, \emph{On existence and uniqueness
  of solutions of {H}amilton-{J}acobi equations}, Nonlinear Anal. \textbf{10}
  (1986), no.~4, 353--370.

\bibitem{LLV1}
Rafael de~la Llave and Enrico Valdinoci, \emph{A generalization of
  {A}ubry-{M}ather theory to partial differential equations and
  pseudo-differential equations}, Ann. Inst. H. Poincar\'e Anal. Non Lin\'eaire
  \textbf{26} (2009), no.~4, 1309--1344.

\bibitem{Ev0}
Lawrence~C. Evans, \emph{On solving certain nonlinear partial differential
  equations by accretive operator methods}, Israel J. Math. \textbf{36} (1980),
  no.~3-4, 225--247.

\bibitem{Ev4}
\bysame, \emph{Periodic homogenisation of certain fully nonlinear partial
  differential equations}, Proc. Roy. Soc. Edinburgh Sect. A \textbf{120}
  (1992), no.~3-4, 245--265.

\bibitem{Hedlund}
Gustav~A. Hedlund, \emph{Geodesics on a two-dimensional riemannian manifold
  with periodic coefficients}, Ann. of Math. \textbf{33} (1932), 719--739.

\bibitem{LPV}
P.-L. Lions, G.~Papanicolau, and S.~R.~S. Varadhan, \emph{Homogenization of
  {H}amilton-{J}acobi equations}, preprint, unpublished, 1987.

\bibitem{MVZ}
Alexandra Monzner, Nicolas Vichery, and Frol Zapolsky, \emph{Partial
  quasimorphisms and quasistates on cotangent bundles, and symplectic
  homogenization}, J. Mod. Dyn. \textbf{6} (2012), no.~2, 205--249.

\bibitem{MoserHD}
J{\"u}rgen Moser, \emph{Minimal solutions of variational problems on a torus},
  Ann. Inst. H. Poincar\'e Anal. Non Lin\'eaire \textbf{3} (1986), no.~3,
  229--272.

\bibitem{Osuna2}
Osvaldo Osuna, \emph{Rigidity of the stable norm on tori}, Rev. Colombiana Mat.
  \textbf{44} (2010), no.~1, 15--21.

\bibitem{ViterboSH}
Claude Viterbo, \emph{Symplectic homogenization}, Preprint arXiv:0801.0206,
  2007.

\end{thebibliography}

\def\cprime{$'$} \def\cprime{$'$} \def\cprime{$'$} \def\cprime{$'$}
\providecommand{\bysame}{\leavevmode\hbox to3em{\hrulefill}\thinspace}
\providecommand{\MR}{\relax\ifhmode\unskip\space\fi MR }
\providecommand{\MRhref}[2]{%
  \href{http://www.ams.org/mathscinet-getitem?mr=#1}{#2}
}
\providecommand{\href}[2]{#2}

\end{document}